\documentclass[acmtoms]{acmtrans2m}

\usepackage{fancyvrb}
\usepackage{algorithm}
\usepackage{amsmath}
\usepackage{amssymb}
\usepackage{psfrag}
\usepackage{graphicx}
\usepackage{color}
\usepackage{url}

\newcommand{\dx}{\, \mathrm{d}x}
\newcommand{\dX}{\, \mathrm{d}X}

\DefineVerbatimEnvironment{code}{Verbatim}{frame=single,rulecolor=\color{blue}}

\begin{document}

\markboth{R. C. Kirby and A. Logg}{Benchmarking domain-specific
  compiler optimizations}

\title{Benchmarking domain-specific compiler optimizations for
  variational forms}
\author{Robert C. Kirby \\
        Texas Tech University
        \and
        Anders Logg \\
        Center for Biomedical Computing, Simula Research Laboratory \\
        Department of Informatics, University of Oslo}

\begin{abstract}
  We examine the effect of using complexity-reducing
  relations~\cite{logg:article:09} to
  generate optimized code for the evaluation of finite element
  variational forms.  The optimizations are implemented in a prototype
  code named FErari, which has been integrated as an optimizing
  backend to the FEniCS Form Compiler,
  FFC~\cite{logg:article:10,logg:article:11}.  In some cases,
  FErari provides very little speedup, while in other cases,
  we obtain reduced local operation counts of a factor of as much as~7.9
  and speedups for the assembly of the global sparse matrix of as
  much as a factor of~2.8 (see Figure~\ref{figure:wadvmat}).
\end{abstract}



\begin{bottomstuff}
Robert C. Kirby, Department of Mathematics and Statistics, Texas Tech University,
Lubbock, TX 79409-1042,
\emph{Email:} \texttt{robert.c.kirby@ttu.edu}.
This work was supported by the
  U.S. Department of Energy Early Career
Principal Investigator Program under award number DE-FG02-04ER25650.
\newline
Anders Logg,
Center for Biomedical Computing,
Simula Research Laboratory,
P.O.Box 134, 1325 Lysaker, Norway.
\emph{Email:} \texttt{logg@simula.no}.
\end{bottomstuff}

\maketitle

\section{Introduction}

Projects such as the FEniCS Form Compiler (hence,
FFC)~\cite{logg:article:10,logg:article:11,logg:www:04},
Sundance~\cite{Lon03,Lon04,www:sundance}, and
deal.II~\cite{www:deal.II} aim to automate important aspects of finite
element computation.  In the case of FFC, low-level code is generated
for the evaluation of element stiffness matrices or their actions,
together with the local-global mapping.  The existence of such a
compiler for variational forms naturally leads one to consider an
\emph{optimizing} compiler for variational forms.  What mathematical
structure in the element-level computations is tedious for humans to
exploit by hand, but possible for a computer to find?  We have
provided partial answers to this question in a series of
papers~\cite{logg:article:07,logg:article:09,KirSco07}. These ideas
have been implemented in a prototype code called FErari, and we
provide an empirical study of the optimizations implemented by FErari
in this paper.  Both FFC and FErari are part of the FEniCS project;
for more information about the software, we refer readers to the
project web page~\cite{www:fenics}.

FFC takes as input a multilinear variational form and generates code
for evaluating that form over affine elements. The formation of the
local stiffness matrix on a single element is expressed as a linear
transformation (known at compile-time) applied to a vector
representing the geometry and coefficient data (known only at
run-time).  The linear transformation depends on the variational form
and finite element basis, but not on the mesh.  This means that the
cost of generating and optimizing the code is independent of the size
of mesh, but depends strongly on the complexity of the variational
form and polynomial degree used.  The generated code is
completely unrolled. This internal kernel is then called for each of
the many elements of the mesh at run-time to compute the global
sparse matrix.
FFC also supports a mode that calls level 2
BLAS~\cite{Dongarra:1988:ESF} rather than generating unrolled code.
This typically gives
comparable run-time performance and smaller executables.  However, the
optimizations we consider here are only possible to apply in the
context of unrolled code.

To a user of FFC, the optimizations are invoked simply with a
\texttt{-O} flag, which turns on a call to FErari and thence a
modified code generator.  It is important to note that the
optimizations considered are similar to, but
typically beyond the abilities of general-purpose compilers to detect.
In assessing the efficacy of these
techniques at reducing run-time, we focus on the construction of the
sparse matrix and its matrix-free application for a variety of
variational forms.  In particular, we study the ``pure'' effect of the
FErari optimizations as well as the optimizations relative to the cost
of inserting into a sparse matrix data structure.

While several fairly theoretical
papers~\cite{logg:article:07,logg:article:09,KirSco07} have shown
that reductions
in arithmetic cost are possible to obtain, there are only very
limited tests of the practical impact of the proposed optimizations.
With some notable exceptions, such as reported in
Figure~\ref{figure:wadvmat} below, the
optimizations provide somewhat disappointing empirical results, such
as only a few percent speedup.  However, it is still important to
include these tests in the literature to bring some completeness to
the theoretical work.  In many cases, the poor speedups are due to
local computation (what we optimize) being dominated by the cost of
insertion into global sparse
data structures.  As memory access is typically very slow compared to
floating point arithmetic, this may not be surprising.
However, it is possible that the optimizations considered here could
perform better in practice in other situations with lower memory
traffic, such as element-by-element or static condensation techniques.
That said, one does obtain significant global speedups in some cases.
For the set of test cases examined below, we obtain a factor of~2.8
global speedup for the assembly of the global sparse matrix of the
weighted advection operator for quartics on tetrahedra (Figure~\ref{figure:wadvmat}).

\section{Finite element assembly and the element tensor}

In finite elements, the nonlinear and linear algebraic problems come
from evaluating the variational forms on the finite element basis
functions.  In our work on FFC and FErari, we have focused on
evaluating multilinear forms over affine elements, and we continue to
do so here.

The typical example is the bilinear form
for Poisson's equation,
\begin{equation}
  a(v, u) = \int_{\Omega} \nabla v \cdot \nabla u \dx.
\end{equation}

If $\{ \phi_j \}_{j=1}^N$ is a finite element basis defined on some
triangulation $\mathcal{T}$ of the domain $\Omega$, the global
stiffness matrix is
\begin{equation}
A_{i} = a(\phi_{i_1},\phi_{i_2}),
\end{equation}
where $i=(i_1,i_2)$ is a multiindex.

The standard algorithm~\cite{ZieTay67,Hug87,Lan99} for computing
the matrix~$A$ is known as \emph{assembly}; it is
computed by iterating over the cells of the mesh
$\mathcal{T}$ and adding from each
cell the local contribution to the global sparse matrix $A$.  A
similar process can compute a global action, in which $A$ is applied
to some vector $u$ without explicitly forming $A$.

The integral defining a multilinear form $a$ may be written as a sum
of integrals over the cells $K$ of a triangulation $\mathcal{T}$ of
the domain $\Omega$:
\begin{equation}
  a = \sum_{K\in\mathcal{T}} a_K,
\end{equation}
and thus
\begin{equation}
  A_{i} = \sum_{K\in\mathcal{T}}
  a_K(\phi_{i_1}, \phi_{i_2}).
\end{equation}
For Poisson's equation, the element bilinear form $a_K$ is thus given
by $a_K(v, u) = \int_{K} \nabla v \cdot \nabla u \dx$.  Finite element
bases are constructed so that each $a_K$ is zero except for a few
basis functions.

For affine elements, as we consider here, the shape functions are
constructed once on a reference element $K_0$ and mapped to each
element of the mesh via an affine mapping $F_K$.  In doing so, one must construct a ``local-global
mapping'' that relates an ordering of the element shape functions to
the global basis functions.  The contribution of element $K$ to the
global matrix $A$ is then evaluated in two stages.  First, a dense element
matrix is computed by evaluating $a_K$ on the shape functions for $K$.
We call this element matrix $A^K$.
Then, each entry of $A^K$ is summed into the appropriate
location in the global sparse matrix as defined by the local-global
mapping.  The first stage is dominated by floating point computation,
the second requires more substantial memory access.

Our work in~\cite{logg:article:10,logg:article:11} has focused on a
general paradigm for efficiently constructing
$A^K$. It has long been known that precomputing certain integrals on the
reference element can speed up computation of the element tensor,
especially for bilinear forms with straight-sided elements.
A general approach to precomputing certain integrals was first introduced in
\cite{KirKne05,logg:article:07} and later formalized and automated in
\cite{logg:article:10,logg:article:11}. A similar
approach was implemented in early versions of
DOLFIN~\cite{logg:preprint:06,logg:www:01,logg:manual:01}, but only
for piecewise linear elements.

 \begin{figure}[htbp]
   \begin{center}
     \psfrag{p0}{$X_1 = (0,0)$}
     \psfrag{p1}{$X_2 = (1,0)$}
     \psfrag{p2}{$X_3 = (0,1)$}
     \psfrag{xi}{$X$}
     \psfrag{x}{$x = F_K(X)$}
     \psfrag{F=}{$F_K(X) = x_1 \Phi_1(X) + x_2 \Phi_2(X) + x_3 \Phi_3(X)$}
     \psfrag{F=}{}
     \psfrag{F}{$F_K$}
     \psfrag{x0}{$x_1$}
     \psfrag{x1}{$x_2$}
     \psfrag{x2}{$x_3$}
     \psfrag{K0}{$K_0$}
     \psfrag{K}{$K$}
     \includegraphics[height=3in]{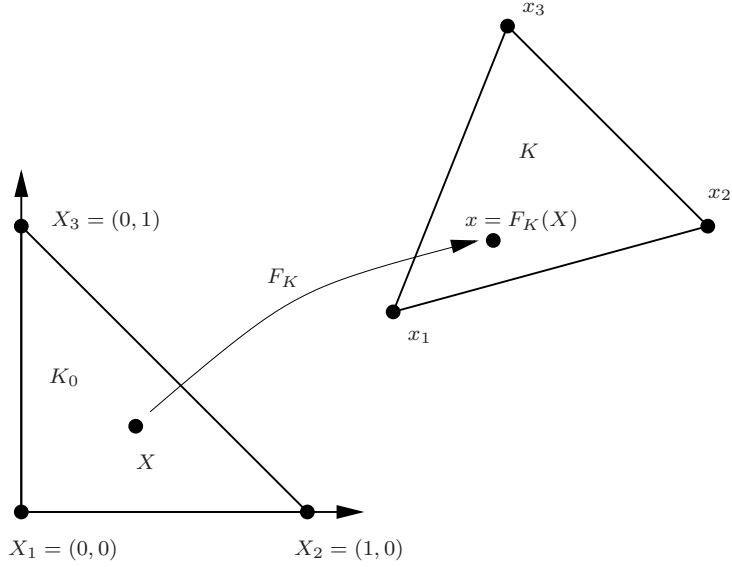}
     \caption{The (affine) mapping $F_K$ from a reference cell $K_0$
       to some cell $K \in \mathcal{T}$.}
     \label{fig:affinemap}
   \end{center}
 \end{figure}

As an example, we consider here the computation of the element matrix
$A^K$ for the Laplacian. When the mapping~$F_K$ from the reference cell
is affine (Figure~\ref{fig:affinemap}), we have for the Laplacian
\begin{equation} \label{eq:poissonA}
  A^K_i =
  \int_K
  \nabla \phi_{i_1}^{K} \cdot
  \nabla \phi_{i_2}^{K} \dx
  =
  \int_K
  \sum_{\beta=1}^d
  \frac{\partial \phi_{i_1}^{K}}{\partial x_{\beta}}
  \frac{\partial \phi_{i_2}^{K}}{\partial x_{\beta}} \dx,
\end{equation}
whence a change of variables yields
\begin{equation}
  A^K_i = \sum_{\alpha\in\mathcal{A}} A^0_{i\alpha} G_K^{\alpha} \quad \forall i \in\mathcal{I}_K,
\end{equation}
where $\mathcal{A}$ and $\mathcal{I}_K$ are sets of allowed multiindices
(depending on the spatial dimension and the discretizing polynomial spaces).
More simply, we can write
\begin{equation} \label{eq:poissoncontraction}
  A^K = A^0 : G_K,
\end{equation}
where
\begin{equation} \label{eq:poissonAandG}
  \begin{split}
    A^0_{i\alpha}
    &=
    \int_{K_0}
    \frac{\partial \Phi_{i_1}}{\partial X_{\alpha_1}}
    \frac{\partial \Phi_{i_2}}{\partial X_{\alpha_2}}
    \dX, \\
    G_K^{\alpha}
    &=
    \det F_K'
    \sum_{\beta=1}^d
    \frac{\partial X_{\alpha_1}}{\partial x_{\beta}}
    \frac{\partial X_{\alpha_2}}{\partial x_{\beta}}.
  \end{split}
\end{equation}
We refer to the tensor~$A^0$ as the \emph{reference tensor}
and to the tensor~$G_K$ as the \emph{geometry tensor}.  For more
details and extensions of this notation to a wide class of multilinear
forms, we refer the reader to our previous
work~\cite{logg:article:10,logg:article:11}.

 In~\cite{logg:article:07,logg:article:09,KirSco07}, we have
explored special mathematical structure that leads to reduced
operation counts.  However, it was studied only in a limited case what
the net impact of FErari optimizations when the cost of global
assembly is counted as well.

\section{A framework for optimization}
\label{sec:matrixvector}

In this section, we present an overview of our framework for
optimization of variational form evaluation. Two different approaches
are presented. The first is a coarse-grained strategy based on
phrasing the tensor contraction~(\ref{eq:poissoncontraction}) as a
matrix-vector or matrix-matrix multiplication that may be computed by
an optimized library call. The second, which is what FErari
implements, exploits the structure of the tensor contraction to find
an optimized computation with a reduced operation count.

\subsection{Tensor contraction as a matrix-vector product}
To evaluate the element tensor~$A^K$, one must evaluate the tensor
contraction~(\ref{eq:poissoncontraction}). A simple approach would be
to iterate over the entries $\{A^K_i\}_{i\in\mathcal{I}_K}$ of $A^K$
and for each entry $A^K_i$ compute the value of the entry by summing
over the set of indices~$\mathcal{A}$.  However, by an
appropriate reshaping of the tensors $A^K$, $A^0$ and $G_K$, one may
phrase the tensor contraction as a matrix--vector product and call an
optimized library routine for the computation of the matrix--vector
product, such as the level 2 BLAS routine DGEMV.  We write
matrix--vector product as
as $a^K = \bar{A}^0 g_K$, where $a^K$ and $g_K$ are $A^K$ and $G_K$
reshaped into vectors and $\bar{A}^0$ is $A^0$ reshaped into a matrix.

Of course, once the computation of one $a^K$ may be computed as a
matrix-vector product, the computation of $\{ a^{K_i} \}_{i=1}^M$ for
some $M$ elements of the mesh can naturally be encoded as a
matrix-matrix multiplication.  Using DGEMM in such a context is an
example of coarse-grained optimization, making good use of cache in a
large computation.  Such an approach necessarily overlooks
problem-specific optimizations such as we find in FErari, but may be
very effective in many circumstances.  It is to be expected that which
approach is preferable will depend strongly on how much structure
FErari finds and how well the resulting algorithms are mapped onto
hardware, as well as whether the computation is large enough for DGEMM
to have good performance.  We do not explore the coarse-grained
strategy further in this paper.

\subsection{Complexity-reducing relations}

The matrix $\bar{A}^0$ is computed at compile-time by FFC, and it
typically possesses significant structure that can be exploited to
reduce the amount of arithmetic needed to multiply it by a vector
$g_K$ at run-time.  It is also helpful to think of the product
$\bar{A}^0 g_K$ as a collection of vector dot products, where vectors
$a^0_i$ are the rows of $\bar{A}^0$.

As an example, we consider forming the weak Laplacian on triangles
using quadratic Lagrange basis functions.  $\bar{A}^0$ is shown in
Table~\ref{tab:bara0}.  We have displayed the index into the
unflattened $A^0$ in the first column, and the rest of row $i$ is the
flattened vector $a_i^0$.  So, the process of forming $A^K$ for some
triangle $K$ is first to compute the geometry vector $g_K$ and then to
form the matrix-vector product $\bar{A}^0 g_K$.  In this case, we will
obtain a vector $a^K$ of length 36, which will be reshaped to the $6\times
6$ element tensor $A^K$.  This is then inserted into the global stiffness
matrix via the local-global mapping.

\begin{table}[htbp]
\begin{center}
\scriptsize
\begin{tabular}{|l|cccc|}
\hline
(0, 0) & 0.5 & 0.5 & 0.5 & 0.5 \\
(0, 1) & 0.16666666667 & 0.0 & 0.16666666667 & 0.0 \\
(0, 2) & 0.0 & 0.16666666667 & 0.0 & 0.16666666667 \\
(0, 3) & 0.0 & 0.0 & 0.0 & 0.0 \\
(0, 4) & 0.0 & -0.66666666667 & 0.0 & -0.66666666667 \\
(0, 5) & -0.66666666667 & 0.0 & -0.66666666667 & 0.0 \\
(1, 0) & 0.16666666667 & 0.16666666667 & 0.0 & 0.0 \\
(1, 1) & 0.5 & 0.0 & 0.0 & 0.0 \\
(1, 2) & 0.0 & -0.16666666667 & 0.0 & 0.0 \\
(1, 3) & 0.0 & 0.66666666667 & 0.0 & 0.0 \\
(1, 4) & 0.0 & 0.0 & 0.0 & 0.0 \\
(1, 5) & -0.66666666667 & -0.66666666667 & 0.0 & 0.0 \\
(2, 0) & 0.0 & 0.0 & 0.16666666667 & 0.16666666667 \\
(2, 1) & 0.0 & 0.0 & -0.16666666667 & 0.0 \\
(2, 2) & 0.0 & 0.0 & 0.0 & 0.5 \\
(2, 3) & 0.0 & 0.0 & 0.66666666667 & 0.0 \\
(2, 4) & 0.0 & 0.0 & -0.66666666667 & -0.66666666667 \\
(2, 5) & 0.0 & 0.0 & 0.0 & 0.0 \\
(3, 0) & 0.0 & 0.0 & 0.0 & 0.0 \\
(3, 1) & 0.0 & 0.0 & 0.66666666667 & 0.0 \\
(3, 2) & 0.0 & 0.66666666667 & 0.0 & 0.0 \\
(3, 3) & 1.3333333333 & 0.66666666667 & 0.66666666667 & 1.3333333333 \\
(3, 4) & -1.3333333333 & -0.66666666667 & -0.66666666667 & 0.0 \\
(3, 5) & 0.0 & -0.66666666667 & -0.66666666667 & -1.3333333333 \\
(4, 0) & 0.0 & 0.0 & -0.66666666667 & -0.66666666667 \\
(4, 1) & 0.0 & 0.0 & 0.0 & 0.0 \\
(4, 2) & 0.0 & -0.66666666667 & 0.0 & -0.66666666667 \\
(4, 3) & -1.3333333333 & -0.66666666667 & -0.66666666667 & 0.0 \\
(4, 4) & 1.3333333333 & 0.66666666667 & 0.66666666667 & 1.3333333333 \\
(4, 5) & 0.0 & 0.66666666667 & 0.66666666667 & 0.0 \\
(5, 0) & -0.66666666667 & -0.66666666667 & 0.0 & 0.0 \\
(5, 1) & -0.66666666667 & 0.0 & -0.66666666667 & 0.0 \\
(5, 2) & 0.0 & 0.0 & 0.0 & 0.0 \\
(5, 3) & 0.0 & -0.66666666667 & -0.66666666667 & -1.3333333333 \\
(5, 4) & 0.0 & 0.66666666667 & 0.66666666667 & 0.0 \\
(5, 5) & 1.3333333333 & 0.66666666667 & 0.66666666667 & 1.3333333333 \\
\hline
\end{tabular}
\caption{The flattened reference tensor for quadratic Lagrange elements on
triangles.  The first column gives the index of the element tensor to
which the row corresponds, and the rest of the columns in the row are the
entries of the flattened vector.}
\label{tab:bara0}
\end{center}
\end{table}

To optimize the evaluation of the element tensor, we look for
dependencies between the vectors $\{a^0_{i}\}_{i\in\mathcal{I}_K}$, or
equivalently the rows of $\bar{A}^0$
that can be used to reduce the cost of forming the matrix-vector
product.  We may only look for structure in
$\{a^0_{i}\}_{i\in\mathcal{I}_K}$, as the $g_K$ vectors are only known
at run-time.  For example, if two
vectors~$a^0_i$ and~$a^0_{i'}$ are collinear (such as the rows (1,0)
and (1,5) in Table~\ref{tab:bara0}), then~$a^0_i \cdot g_K$
may be computed using~$a^0_{i'}\cdot g_K$ in only one multiply, and
vice versa.  If the Hamming distance (number of different entries
between~$a_0^i$ and~$a_0^{i'}$) is $k$, then the result~$a^0_{i'}\cdot
g_K$ can be computed from~$a^0_{i}\cdot g_K$ in about $k$ multiply-add
pairs, and vice versa.  These kinds of relations are called
``complexity-reducing relations'', and they are related to common
subexpressions.  Note that using such a relationship requires that the
code for the dot products be unrolled.  As with FFC, there may come a
point at which code bloat outweighs gains in arithmetic cost, but we
remark that code optimized by FErari contains fewer arithmetic
operations and hence is smaller than the standard FFC output, but much
larger than using the BLAS mode of FFC.

In~\cite{logg:article:09}, we constructed a weighted, undirected
graph, the vertices of which were the vectors~$a^0_{i}$ and the
weights of whose edges were the pairwise distances under a
complexity-reducing relation (the cost of computing one entry in the
element matrix from another).  We proved that a minimum spanning tree
of this graph encodes a minimal-arithmetic (in a specific sense)
algorithm for evaluating the product of~$\bar{A}^0$ with an arbitrary
input vector.

\begin{figure}
\begin{center}
\includegraphics[height=7in]{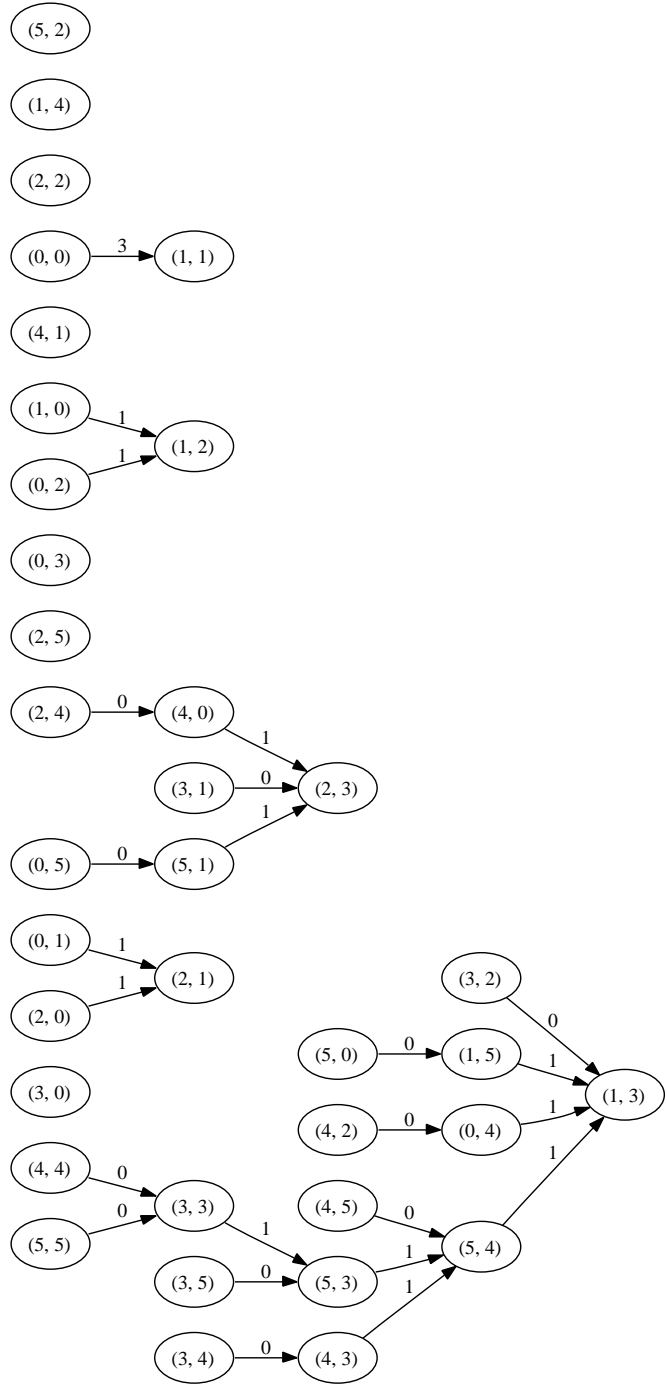}
\caption{Dependency graph for forming the element stiffness matrix for the
Laplacian using quadratic Lagrange triangles as
determined by FErari.}
\label{fig:msf}
\end{center}
\end{figure}

In Figure~\ref{fig:msf}, we show the dependency graph generated by
FErari. The arrows indicate dependency rather than implication.  That
is, the arrow from (0,0) to (1,1) indicates that the result of
computing $a^0_{(1,1)} g_K$ is used to compute $a^0_{(0,0)} g_K$.
Hence, the implied flow of computation is from right to left, and
disconnected components in the graph are independent of each other.

As one extension of this technique, we notice that many of the vectors
may be computed effectively by ignoring multiplication by zero.  For
example, entry (1,3) in Table~\ref{tab:bara0} only has one nonzero
entry.  It makes sense to generate code for forming $a^0_{(1,3)} g_K$
explicitly instead of using a complexity-reducing relation. In this case,
we have ``snipped'' the edge from the entry (1,3) to its parent in the
minimum spanning tree before generating code and thus this entry has
no outgoing arrows. Hence, we properly have a forest rather than a
tree.

Many other kinds of structure may be found in $\bar{A}^0$.  For
example, in many cases one can prove that the $g_K$ tensor has
symmetries along certain axes.  We used this, for example,
in~\cite{logg:article:07,logg:article:09}, but have yet
not automated the detection of such structure.  Also, frequently
three or more rows of $\bar{A}^0$ will be linearly dependent.  A first
attempt at exploiting this structure is found in~\cite{KirSco07}, but
our present work is limited to complexity-reducing relations.

\section{Benchmark results}
\label{sec:benchmark}

For a range of forms and polynomial degrees, we report several
quantities for forming the matrix and its action. First, we report
the base operation count $|\mathcal{I}_K| \, |\mathcal{A}|$ for
forming the element tensor $A^K$, as well as the operation counts generated by
FFC\footnote{FFC reduces the base operation count by omitting
computation of zeros when the element tensor is sparse.} and the
FErari optimizations.  Having generated code for the local element
computation from both FFC and FErari, we compare the run-time for
these codes being executed several times.  This measures the efficacy
of FErari at exactly the point it seeks to optimize.  Then, to provide
a broader context, we present the speedup obtained in the global
assembly process, when the overhead of sparse data structures is
included.

In each case, we generated code for the local and global computation
both with and without FErari optimizations. This code was compiled and
run on an IBM Thinkpad T60p with 2GB of RAM and a dual core Intel
T2600 chip running at 2.16 GHz. The operating system was Ubuntu Linux
with kernel 2.6.17-10-386.  The compiler was \texttt{g++} version
4.1.2 using optimization flag \texttt{-O2} on all variational forms
except the weighted Laplacian operator and action using quartics in
3D.  The compiler and machine could only handle optimization mode
\texttt{-O0} in these cases. This illustrates a challenge with our
approach to finite element code generation based on the tensor
representation~(\ref{eq:poissoncontraction}). Since
straight-line code is generated for the computation of the element
tensor, complicated forms or high-dimensional finite element spaces
may lead to generation of large amounts of code which the C++ compiler
is not able to handle, particularly in optimized mode. For these
forms, generating code based on quadrature rather than tensor
contraction with FFC/FErari could be more practical.

For two-dimensional problems, we used a
regular triangulation based on subdividing a $64 \times 64$ square
mesh into right triangles, resulting in a total of 4,225 vertices and
8,192 triangles.  For three dimensions, we used a $16 \times 16 \times
16$ partition of the unit cube into 4,913 vertices and 24,576
tetrahedra. The timing was performed adaptively to ensure that at
least one second of CPU time elapsed for a set of at least ten
repetitions for each test case.  For the sparse matrix data structure,
a simple \texttt{std::vector<std::map<unsigned int, double> >} was
used, which was found competitive with insertion into a sparse PETSc
matrix.

In most cases, we find decent speedup in the operation count, although
it does not always translate into a speedup in the runtime for the
local computation.  FErari is currently architecture-unaware.
Rearranging the matrix-vector computation in a way that makes poor use
of registers, for example, can more than offset reductions in the actual
amount of arithmetic.  A better result would be obtained by
somehow combining the graph-based optimizations with an architecture
model, or using a special-purpose compiler such as Spiral~\cite{Pueschel:05}.

Moreover, even a speedup in local computation does not always improve
the global cost of assembling a matrix or vector.  If a relatively
small amount of work is required to compute $A^K$, then the cost of
assembling it into the global matrix or vector may dominate;
reductions in arithmetic are not significant.  On the other hand, when
the construction of $A^K$ is relatively expensive, then speedup in the
construction of the global matrix or vector can be realized by
reduction of arithmetic in the local computation. In our empirical
results, we observe a tendency of FErari to provide better global
speedups for more complicated variational forms.

\subsection{Laplacian}

First, we consider the Laplacian, with the variational form
\begin{equation}
  \label{eq:laplacian}
  a(v, u) = \int_{\Omega} \nabla v \cdot \nabla u \dx.
\end{equation}
We use Lagrange polynomials $P_k$ of degree $k = 1, 2, \ldots, 5$
on triangles and degree $k = 1, 2, \ldots, 4$ on
tetrahedra.\footnote{The polynomial degree on tetrahedra was limited
by available resources to compute the optimization.}

\begin{figure}
  \centerline{\includegraphics[height=3in]{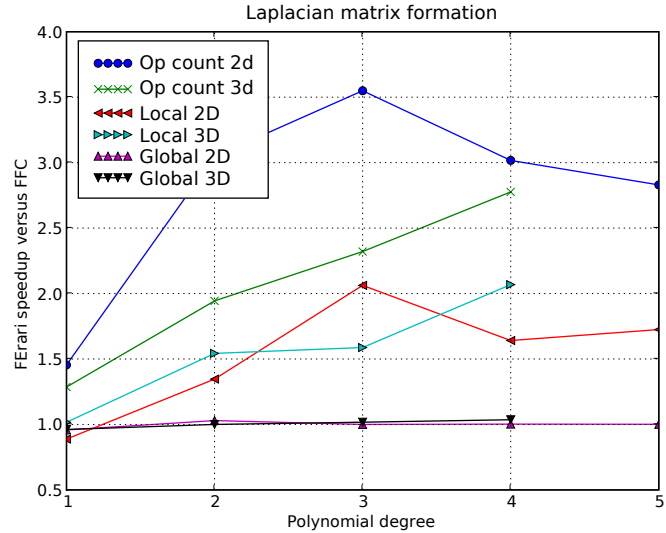}}
  \caption{Speedup in operation count, local run-time and global run-time for using FErari versus FFC only
    for the Laplacian~(\ref{eq:laplacian}).}
  \label{figure:lapmat}
\end{figure}

 In each case, FErari provides up to about a factor of three
improvement in operation count. The reduction in operation count,
local computation time, and global computation time required is
plotted in Figure~\ref{figure:lapmat}.  The reduction in arithmetic
reduces the run-time to evaluate the local stiffness matrix
(multiplying by~$\bar{g}_K$) by a factor of 1.5 to 2 in both two and
three dimensions. However, the reduction does not have a major impact
on the global time to assemble the matrix.  In this case, there are
very few arithmetic operations needed to construct the local matrix,
and the cost of inserting into the global matrix overshadows the gains
FErari provides.

\begin{figure}
  \centerline{\includegraphics[height=3in]{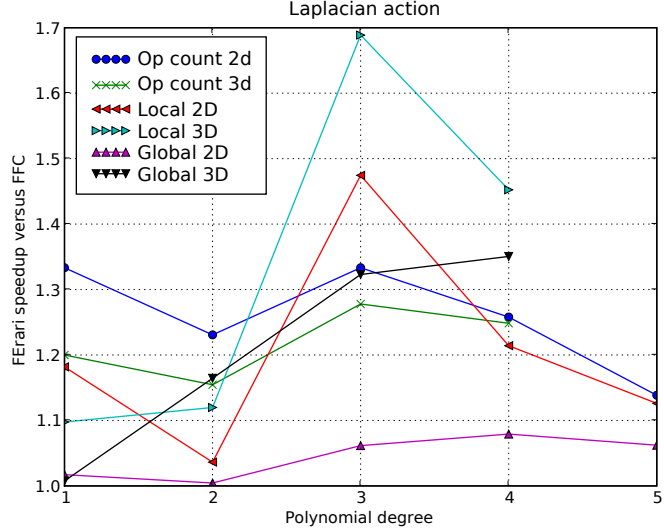}}
  \caption{Speedup in operation count, local run-time and global run-time for using FErari versus FFC only
    for the action of the Laplacian~(\ref{eq:laplacian}).}
  \label{figure:lapact}
\end{figure}

We also consider the matrix action as needed in a Krylov
solver. Assembling into a global vector is less expensive than into a
global matrix, and we see better speedups in evaluating the action of
the Laplacian operator. In this case, FFC and FErari generate code for
evaluating~(\ref{eq:laplacian}) with $u$ a member of the finite
element space.  Speedup of this operation is felt at each iteration of
a Krylov method and so translates directly into decreased solve time.
The matrix $\bar{A}^0$ has the same entries as for forming
the stiffness matrix, but has a different shape.  In this case, the
shape is $|P_k| \times (d^2 |P_k|)$.
Note that FErari does not do as well for the action as for forming the
matrix.  Although the entries of $\bar{A}^0$ are the same as before,
the difference in shapes complicates finding collinear relationships.
When the rows have only $d^2$ (4 or 9) entries for the stiffness
matrix, more collinearity is found than when there are $|P_k|$ times as
many entries.  However, finding Hamming distance relations is as
effective as before.  Despite the smaller reduction in operation
count, the effect of the optimizations on run-time is much greater
than in forming the matrix, as we can see by comparing
Figure~\ref{figure:lapact} to Figure~\ref{figure:lapmat}. A global
speedup of about 10\% is observed for degrees three through five in
two dimensions, and a speedup of 20\%--40\% for quadratics through
quartics in three dimensions.  Again, only a small improvement is
observed for low order methods.

\subsection{Weighted Laplacian}

Now, we consider the form
\begin{equation}
  \label{eq:laplacian,weighted}
  a(v, u, w) = \int_{\Omega} w \nabla v \cdot \nabla u \dx,
\end{equation}
for a fixed weight $w$ where we assume that $v,u,w$ all come from the
same Lagrange finite element space. In this case, the presence of the
coefficient $w$ makes the local form more expensive to evaluate.  The
matrix $\bar{A}_0$ now has $|P_k|^2$ rows and $d^2 |P_k|$
columns. However, the graph of the global matrix for this form is the
same as for the constant coefficient case, assuming the same basis and
mesh are used.  Consequently, the cost of assembly is exactly the same
once $A^K$ is constructed.

\begin{figure}
  \centerline{\includegraphics[height=3in]{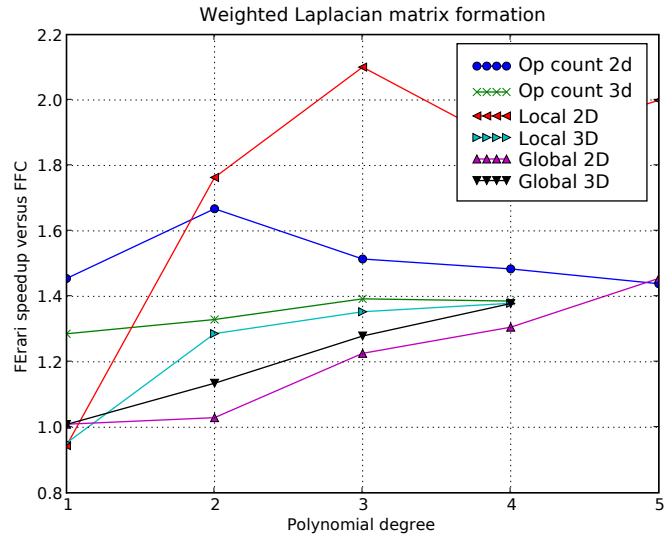}}
  \caption{Speedup in operation count, local run-time and global run-time for using FErari versus FFC only
    for the weighted Laplacian~(\ref{eq:laplacian,weighted}).}
  \label{figure:wlapmat}
\end{figure}

Again, FErari reduces
the operation count and run-time for the local computation considerably.
Given that the arithmetic cost is much larger than for the
constant-coefficient case, it is not surprising that the global
speedups are much better, as seen in Figure~\ref{figure:wlapmat}.

\begin{figure}
  \centerline{\includegraphics[height=3in]{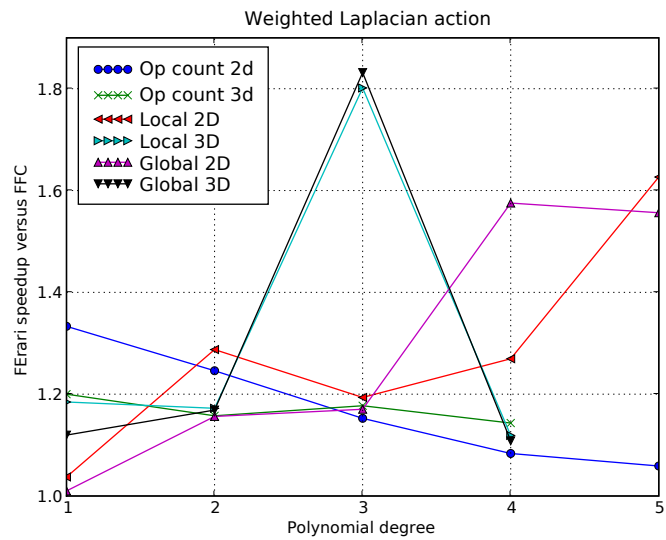}}
  \caption{Speedup in operation count, local run-time and global run-time for using FErari versus FFC only
    for the action of the weighted Laplacian~(\ref{eq:laplacian,weighted}).}
  \label{figure:wlapact}
\end{figure}

As before, $\bar{A}^0$ has the same entries but a different shape when
the action of the form is considered. Now, the shape is $|P_k| \times
(d^2 |P_k|^2)$.
While FErari does not reduce the operation
count for the matrix action as significantly as it does for the matrix
itself, the global speedups are more significant (Figure~\ref{figure:wlapact}).

\subsection{Advection}

Next, we consider the advection operator
\begin{equation}
  \label{eq:adv}
  a(v, u) = \int_{\Omega} v (\beta \cdot \nabla u) \dx,
\end{equation}
where~$\beta$ is some constant vector and consider forming the global
stiffness matrix and its action.  For the matrix, the dimension of
$\bar{A}^0$ is $|P_k|^2 \times d^3$.
The advection $\beta$ is defined as a
piecewise constant vector-valued Lagrange function which has $d$
degrees of freedom on each element. As a result, the matrix
$\bar{A}^0$ is physically of dimension $|P_k|^2 \times d^3$, but the
number of nonzero elements scales like $|P_k|^2 \times d^2$. This is
because the reference tensor $A^0$ generating the matrix $\bar{A}^0$
is formed as an outer product with $\Phi_{\alpha_1}[\alpha_2] =
\delta_{\alpha_1\alpha_2}$, that is, component $\alpha_2$ of the
piecewise constant vector-valued basis function $\Phi_{\alpha_1}$.
Precontracting the reference tensor along dimensions $\alpha_1,
\alpha_2$ would thus reduce the size of the matrix $\bar{A}^0$ to $|P_k|^2
\times d^2$. Low-order elements like piecewise constants and linears
often generate particular structures that can be used for further
optimizations. Such optimizations are not handled by FErari and are an
interesting venue for further research.

\begin{figure}
  \centerline{\includegraphics[height=3in]{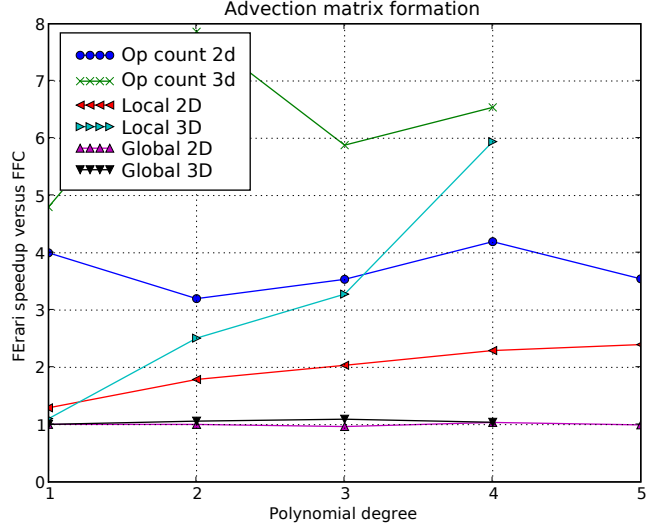}}
  \caption{Speedup in operation count, local run-time and global run-time for using FErari versus FFC only
    for the advection operator~(\ref{eq:adv}).}
  \label{figure:advmat}
\end{figure}

\begin{figure}
  \centerline{\includegraphics[height=3in]{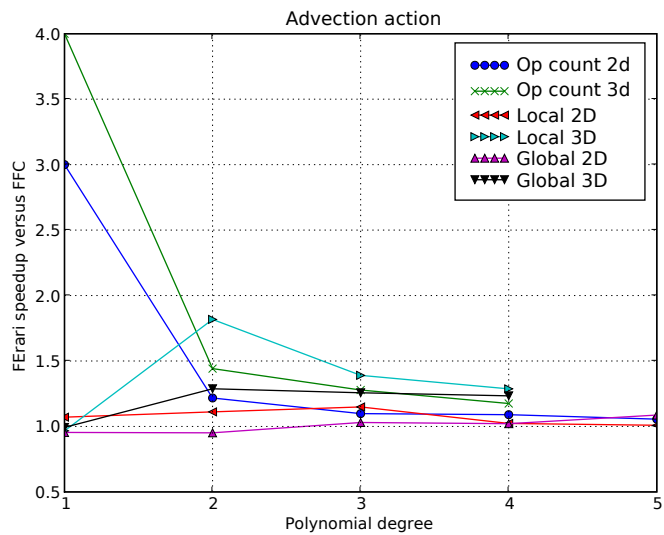}}
  \caption{Speedup in operation count, local run-time and global run-time for using FErari versus FFC only
    for the action of the advection operator~(\ref{eq:adv}).}
  \label{figure:advact}
\end{figure}

As with forming the Laplacian, the reduced operation counts
do not significantly affect the global
runtime (Figure~\ref{figure:advmat}). The operation
counts and speedups for the matrix action are found in
Figure~\ref{figure:advact}. Global
speedup is again most significant for higher order elements in three
dimensions.

\subsection{Weighted advection in a coordinate direction}

\begin{figure}
  \centerline{\includegraphics[height=3in]{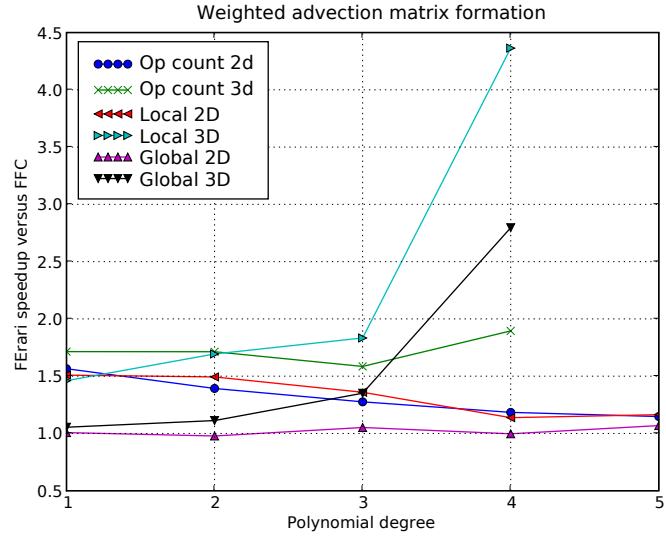}}
  \caption{Speedup in operation count, local run-time and global run-time for using FErari versus FFC only
    for the weighted advection operator~(\ref{eq:wadv}).}
  \label{figure:wadvmat}
\end{figure}

\begin{figure}
  \centerline{\includegraphics[height=3in]{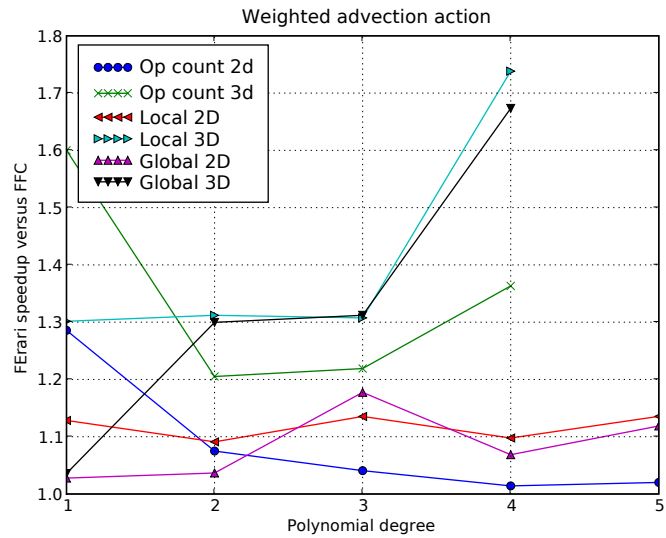}}
  \caption{Speedup in operation count, local run-time and global run-time for using FErari versus FFC only
    for the action of the weighted advection operator~(\ref{eq:wadv}).}
  \label{figure:wadvact}
\end{figure}

Finally, we consider the advection operator oriented along a coordinate
axis, but with the velocity field varying in space (projected into the
finite element space):
\begin{equation}
  \label{eq:wadv}
  a(v, u, w) = \int_{\Omega} v w \frac{\partial u}{\partial x_1} \dx,
\end{equation}
We consider forming the matrix and its action for a fixed weight
$w$. This operator is a portion of the trilinear momentum advection
term in the Navier--Stokes equations. For constructing the matrix, we
observe a nice speedup in local computation, although in two
dimensions this has only a marginal effect on the global run-time for
assembly.  However, we gain significantly for higher-order elements in
three dimensions, where we see a global speedup with 180\% (a factor
2.8) for quartics. The operation counts for the local matrix
construction and action are shown
in Figures~\ref{figure:wadvmat} and~\ref{figure:wadvact}.

\subsection{Speedup versus work}

\begin{figure}
  \centerline{\includegraphics[height=3in]{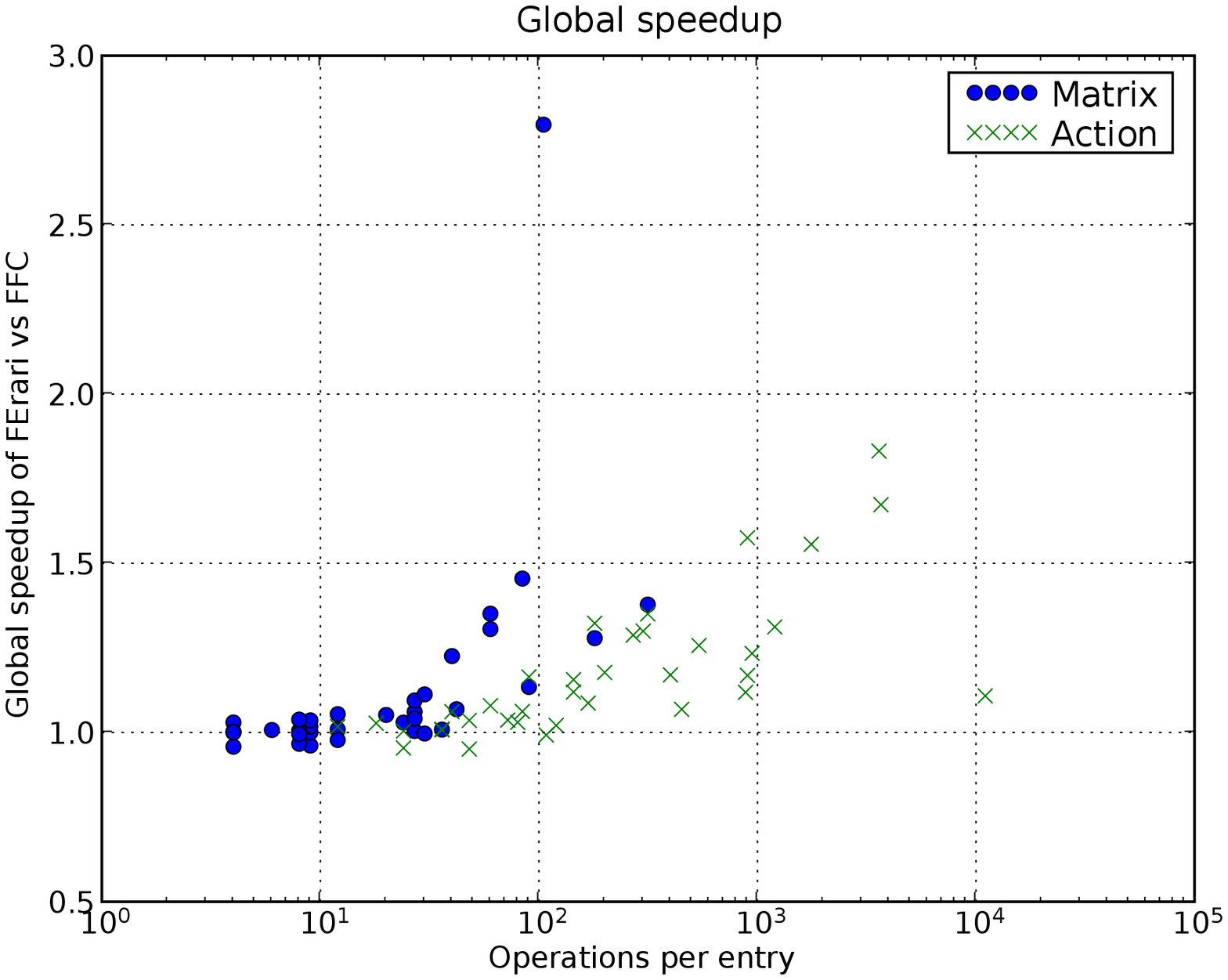}}
  \caption{The global speedup that FErari produces over FFC is plotted against
    the number of columns in the associated reference
    matrix $\bar{A^0}$, which is a measure of the work required to
    compute each entry of $A^K$.}
\label{fig:speedup}
\end{figure}

As we noted before, reducing floating-point arithmetic is expected to
be more significant to the global computation when the individual
entries in the local matrix or vector are already expensive to
compute.  As a test of this, we plot the speedup of FErari over FFC
against the number of columns in each reference operator $\bar{A}^0$
in Figure~\ref{fig:speedup}. We do this for all orders and forms,
considering matrices and their actions separately.  Although it is not
an exact relation (as to be expected), Figure~\ref{fig:speedup}
does indicate a general trend of speedup increasing with the base cost
of work per entry.

\subsection{Compile times}

It is important to quantify the additional compile-time cost of using
FErari within FFC.  In some situations, especially in a just-in-time
compilation, the significant additional cost will outweigh the
potential run-time gains.  In this section, we report compile times for
a few forms as an example.  It should be remembered, however, that
FErari is currently implemented in Python and far from tuned for
performance. A better implementation should improve these compile
times.

Tables~\ref{tab:compile,nonoptimized} and~\ref{tab:compile,optimized}
give the compile times for FFC without and with FErari optimizations
respectively. We also report the time for compiling the C++ code
generated by FFC with GCC (\texttt{g++}). We note a few interesting
details from these numbers. First, we note that the FErari
optimizations may take considerable time, in particular for high
degree polynomials and forms containing coefficients. Further, we note
that it may also take considerable time to compile the generated
code. Finally, we note that GCC may in some cases run faster if the
generated code has already been optimized by FErari.  This gain is
small compared to the cost of running FErari, and is directly
attributable to the resulting unrolled code having fewer operations.

\begin{table}
  \begin{center}
    \footnotesize
    \begin{tabular}{|lc|c|c|c|}
      \hline
      Form & Degree & FFC & GCC & GCC -O2 \\
      \hline
      Laplacian operator & 1 & 0.016 & 2.3 & 2.3 \\
      Laplacian operator & 2 & 0.035 & 2.2 & 2.5 \\
      Laplacian operator & 3 & 0.13 & 2.5 & 3.7 \\
      Weighted Laplacian operator & 1 & 0.029 & 2.2 & 2.4 \\
      Weighted Laplacian operator & 2 & 0.26 & 2.8 & 5.2 \\
      Weighted Laplacian operator & 3 & 2.3 & 9.1 & 130 \\
      \hline
    \end{tabular}
    \caption{Compile times in seconds for FFC, GCC and GCC with optimization \texttt{-O2} for
      a set of forms.}
    \label{tab:compile,nonoptimized}
  \end{center}
\end{table}

\begin{table}
  \begin{center}
    \footnotesize
    \begin{tabular}{|lc|c|c|c|}
      \hline
      Form & Degree & FFC -O & GCC & GCC -O2 \\
      \hline
      Laplacian operator & 1 & 0.12 & 2.1 & 2.3 \\
      Laplacian operator & 2 & 4 & 2.2 & 2.5 \\
      Laplacian operator & 3 & 68 & 2.4 & 3.3 \\
      Weighted Laplacian operator & 1 & 0.23 & 2.2 & 2.4 \\
      Weighted Laplacian operator & 2 & 22 & 2.6 & 4.5 \\
      Weighted Laplacian operator & 3 & 760 & 7.2 & 78 \\
      \hline
    \end{tabular}
    \caption{Compile times in seconds for FErari-optimized FFC, GCC and GCC with optimization \texttt{-O2} for
      a set of forms.}
    \label{tab:compile,optimized}
  \end{center}
\end{table}

\section{Conclusions}

Several things emerge from our empirical study of optimizing FFC with
FErari.  In certain contexts, FErari can provide tens of percent
to a few times speedup in runtime in forming or applying stiffness
matrices.  Moreover, these cases tend to be the computationally harder
ones (three dimensions, higher order polynomials). However, FErari is
not without its costs.  It dramatically adds to the compile-time for
FFC, and when used for simple forms can actually hinder runtime.

Besides improving the run-time performance of finite element codes
generated by FFC and FErari, our results shed some light on where
FErari could be improved and in how a fully functional optimizing
compiler for finite elements might be developed.  First, our
calculations did little to optimally order the degrees of freedom;
better ordering algorithms should decrease the cost of insertion.
Second, algorithms trying to maximize performance must have some
awareness of the underlying computer architecture.  The success of
Spiral in signal processing suggests this should be possible.
Moreover, knowing when to do what kinds of optimization, such as
FErari's fine-grained optimization versus a coarse-grained level~3
BLAS approach, must be determined.  This must also be compared against
when quadrature-based algorithms might be effective, as well as
whether the stiffness matrix should be explicitly constructed,
statically condensed, or applied without being constructed.

\subsection*{Acknowledgments}

Logg is supported by an Outstanding Young Investigator grant
from the Research Council of Norway, NFR 180450.

\bibliographystyle{acmtrans}
\bibliography{bibliography}

\end{document}